\documentclass{didier}

\title{Sums of even ascending powers}
\author{Chan Ieong Kuan, Didier Lesesvre, Xuanxuan Xiao}
\date{\today}

\address{School of mathematics (Zhuhai) \newline
Zhuhai Campus, Sun Yat-Sen University \newline
Tangjiawan, Zhuhai, Guangdong, 519082, China (PRC)}
\email{kuanchi3@mail.sysu.edu.cn}

\address{School of mathematics (Zhuhai) \newline
Zhuhai Campus, Sun Yat-Sen University \newline
Tangjiawan, Zhuhai, Guangdong, 519082, China (PRC)}
\email{lesesvre@math.cnrs.fr}

\address{Macau University of Science and Technology \newline
Faculty of Information Technology \newline
Taipa, Macau}
\email{xxxiao@must.edu.mo}

\usepackage{stmaryrd}

\newcommand{\N}{\mathbf{N}}

\begin{document}

\maketitle

\begin{abstract}
Freiman and Scourfield proved that any large enough integer can be written as a sum of a certain number of ascending even powers. We use the circle method to provide the first explicit bound on this number, and show that any large enough integer can be written as a sum of 133 ascending even powers. 
\end{abstract}

\setlength{\parskip}{0em}
\setcounter{tocdepth}{1}
\setlength{\parskip}{0.5em}

\section{Introduction}
\label{sec:introduction}

\subsection{Ascending powers}

Since the first steps of the circle method at the beginning of the twentieth century, a wide literature has been published concerning the Waring problem and its generalizations of mixed types. The question is the one regarding the possibility of representing a large enough integer $n$ in the form 
\begin{equation}
\label{ascending-powers-problem}
n = x_1^{k_1} + \cdots + x_s^{k_s}, \qquad k_1 \leqslant \cdots \leqslant k_s, \qquad x_i, \, k_i \in \N.
\end{equation}

Freiman \cite{freiman_solution_1949} and Scourfield \cite{scourfield_generalization_1960} characterized the sequences of powers for which \eqref{ascending-powers-problem} is solvable for all large enough integer $n$. More precisely, they established the following result.

\begin{thm}[Freiman-Scourfield]
Let $(k_i)_i$ be a non-decreasing sequence of positive integers. The series $\sum k_i^{-1}$ is divergent if and only if for any $j \in \N$, there exists an $s \geqslant j$ such that every large enough integer $n$ is representable as
\begin{equation}
\label{freiman-scourfield}
n = x_j^{k_j} + x_{j+1}^{k_{j+1}} +\cdots + x_s^{k_s}, \qquad x_j, \ldots, x_s \in \N.
\end{equation}
\end{thm}

Despite this result, few is known concerning the lowest possible $s$ for which \eqref{freiman-scourfield} holds provided the series $\sum k_i^{-1}$ is divergent.  Many authors, among which Roth \cite{roth_proof_1949, roth_problem_1951}, Thanigasalam \cite{thanigasalam_additive_1968, thanigasalam_sums_1980}, Vaughan \cite{vaughan_warings_1991, vaughan_further_1995}, \linebreak Brüdern  \cite{brudern_sums_1987, brudern_problem_1988} and Ford \cite{ford_representation_1995, ford_representation_1996}, struggled for half a century to determine the least $s$ for which any large enough integer can be written as a sum of ascending powers,
\begin{equation}
\label{full-ascending-powers}
n = x_2^2 + x_3^3 + \cdots + x_s^s, \qquad x_2, \ldots, x_s \in \mathbf{N}.
\end{equation}

The purpose of this paper is to address another variation of the Waring problem of mixed type, restricted to the sequence of ascending even powers $\{ 2n \}_{n \in \N}$ . This elegant variation has been recently explored by Brüdern. The set of integers representable as 
\begin{equation}
n = x_2^2 + x_4^4 + x_6^6, \qquad x_2, x_4, x_6 \in \N,
\end{equation}

is of zero density since the sum of reciprocal exponents fails to reach one. The main result of Brüdern's paper \cite{brudern} implies in particular that the set of integers representable as 
\begin{equation}
n = x_2^2 + x_4^4 + x_6^6 + x_8^8, \qquad x_2, x_4, x_6, x_8 \in \N,
\end{equation}

has positive density, but it does not reach 1. It turns out that a corollary of his result establishes that the set of integers representable as 
\begin{equation}
n = x_2^2 + x_4^4 + x_6^6 + x_8^8 + x_{10}^{10}, \qquad x_2, x_4, x_6, x_8, x_{10} \in \N,
\end{equation}

is of density one. This settles the problem of the even ascending powers in the density aspect. The purpose of this paper is to give a bound on the number of even ascending powers necessary to write all but finitely many natural numbers in this form. Investigating this variation has shed some light on interesting features concerning the choices of parameters in the circle method, which are often kept undisclosed, and we hope this paper can be useful in this pedagogical sense. We put a particular emphasis on clarifying these choices and commenting the origins of each parameter and on the limitations for the final value of $s$. These heuristics facts are labeled as remarks all along the paper. The main result is the following. 

\begin{thm}
\label{thm:even-ascending-powers}
Let $s=133$. Every sufficiently large natural number $n$ is representable in the form
\begin{equation}
\label{even-ascending-powers}
n = {x_2}^2 + {x_4}^4 + {x_6}^6 + \cdots + x^{2s}_{2s},  \qquad x_2, \ldots, x_{2s} \in \N.
\end{equation}
\end{thm}

Let $R_s(n)$ be the number of ways of writing $n$ as in \eqref{even-ascending-powers}. The aim is to prove that $R_s(n) > 0$ for $n$ large enough, and to find the best possible $s$ for which it happens.

\textit{Remark.} The theoretical limit of the circle method is given by a sum of reciprocal exponents equal to 2. The quality of the result can henceforth be judged by how close to 2 is the sums of reciprocal exponents. In the case of the growing powers, Ford \cite{ford_representation_1996} reaches the value
\begin{equation}
\label{sum-reciprocals-ford}
\sum_{k=2}^{15} \frac{1}{k} \approx 2.21,
\end{equation}

\noindent while the result presented here in the case of the even ascending powers yields
\begin{equation}
\mu = \sum_{k=1}^{133} \frac{1}{2k} \approx 2.73.
\end{equation}

This is slightly worse than Ford's case, and can be understood by the fact that our sequence is growing faster than the one in \eqref{sum-reciprocals-ford} and some specific bounds for small powers providing strong savings in Ford's case cannot be used in the present case. We prove in fact a quantitative version of Theorem \ref{thm:even-ascending-powers}.

\begin{thm}
\label{thm:3}
Let $s = 133$. We have
\begin{equation}
R_s(n) \gg F(0) n^{-1} \asymp n^{\mu - 1}.
\end{equation}
\end{thm}

\subsection{Outlook of the proof}

This is a typical problem of additive number theory. Similar to the Waring problem or to the full ascending power problem \eqref{full-ascending-powers}, our problem is amenable to the Hardy-Littlewood circle method. For the most part, the numerous parameters are kept general until explicit numerical choices are needed for conclusions, hence motivating the choices made all along. We define the Farey dissection and the major arcs in Section \ref{sec:circle-method}, estimate the contribution of minor arcs in Section \ref{sec:minor-arcs}, approximate the generating functions on major arcs in Section \ref{sec:approximation}, prune the major arcs to logarithmically-scaled arcs in Section \ref{sec:pruning} and finally conclude by estimating the main term in Section \ref{sec:singular}.

\subsection{Acknowledgments}

The first-named author is supported in part by NSFC (No.11901585). The second-named author is infinitely indebted to J\"{o}rg Br\"{u}dern for having been a strong source of motivation all along the project as well as for enlightening discussions. The thrid-named author is supported in part by NSFC (No. 11701596) and the Science and Technology Development Fund, Macau SAR (No. 0095/2018/A3).

\section{Circle method}
\label{sec:circle-method}

\subsection{Definition of major arcs}

Let $n \geqslant 1$ and $0 < \tau < 1/2$. Set $X=n^\tau$. Denote $\|x\|$ the integer norm of $x$, that is the distance to the closest integer. For $1 \leqslant a \leqslant q \leqslant X$ with $(a,q) = 1$, introduce the major arc around $a/q$ defined by
\begin{equation}
\mathfrak{M}(X; q, a) = \left\{ \alpha \in [0,1] \ : \left\| \alpha - \frac{a}{q} \right\| \leqslant \frac{X}{qn} \right\}.
\end{equation} 

The major arcs $\mathfrak{M}(X)$ are the union of all these arcs for $q \leqslant X$, namely
\begin{equation}
\mathfrak{M}(X) = \bigcup_{q \leqslant X} \bigcup_{\substack{1 \leqslant a \leqslant q \\ (a,q) = 1}} \mathfrak{M}(X ; q, a).
\end{equation}

The major arcs $\mathfrak{M}(X; q, a)$ are pairwise disjoint provided $X < \frac{1}{2} n^{1/2}$, so that any $\alpha \in \mathfrak{M}(X)$ uniquely determines the associated $a$ and $q$. We use this fact every time these quantities appear without notice and $\alpha$ is fixed in a major arc. Define the minor arcs as the remaining points on the circle, that is to say
\begin{equation}
\mathfrak{m}(X) = [0,1] \backslash \mathfrak{M}(X).
\end{equation}

\subsection{Analytic generating functions}

Introduce the dyadic exponential sums, for all $k \geqslant 1$ and $\alpha \in \R$,
\begin{equation}
f_k(\alpha) = \sum_{m \in X_k} e\left(\alpha m^k\right), \qquad \text{where } \quad X_k = \left\llbracket n^{1/k} + 1, 2n^{1/k}  \right\rrbracket,
\end{equation}

\noindent and their smooth analogue, for a certain $\gamma > 0$ to be determined later,
\begin{equation}
g_k(\alpha) = \sum_{m \in Y_k} e\left(\alpha m^k\right), \qquad \text{where } \quad Y_k = \mathcal{A}\left(n^{1/k}, n^\gamma\right).
\end{equation}

Here, $\mathcal{A}(n^{1/k}, n^\gamma)$ stands for the $n^\gamma$-smooth numbers less than $n^{1/k}$, that is to say
\begin{equation}
\mathcal{A}\left(n^{1/k}, n^\gamma\right) = \left\{ x \leqslant n^{1/k} \ : \ p | x \Longrightarrow p \leqslant n^\gamma \right\}.
\end{equation}

Let $K=\{2, 4, \ldots, 2s\}$ be the set of indices we consider in \eqref{even-ascending-powers}. For the partition $K = K_f \sqcup K_g$, introduce the generating function
\begin{equation}
F = \prod_{k \in K_f} f_k \prod_{k \in K_g} g_k.
\end{equation}

In particular, the associated Fourier coefficients
\begin{equation}
\label{rN}
r_s(n) = \int_{\R/\Z} F(\alpha)e(-n\alpha) \mathrm{d}\alpha, \qquad n \in \mathbf{N},
\end{equation}

\noindent are less than $R_s(n)$, the number of solutions to the initial problem \eqref{even-ascending-powers}. It is therefore sufficient to prove that $r_s(n) > 0$ for $n$ large enough in order to prove Theorem \ref{thm:even-ascending-powers}. The expected size is of order $F(0)n^{-1}$, as proven in the last section and stated in Theorem \ref{thm:3}. This provides a guide to estimating error terms.

\section{Contribution of minor arcs}
\label{sec:minor-arcs}

\subsection{Strategy}

We begin by bounding the contribution of minor arcs to the integral \eqref{rN}. Let $K_f = \{2, 4, 6\}$ and $K_g = \{8, \ldots, 2s\}$. Introduce the set of indices
\begin{align*}
K & = \{2\} \sqcup K_1 \sqcup K_2.
\end{align*}

\textit{Remark.} We would ultimately like to take $K_f = \emptyset$, however non-smooth functions $f_k$ are necessary in order to apply the iterative methods. What determines the choice of $K_1$ and $K_2$ is an optimization between both sets in order to get the best values in the application of the mixed mean value algorithm. See Sections \ref{subsec:minor} and \ref{f4} for more details.

Introduce
\begin{align*}
F_1 & =  \prod_{k \in K_1 \cap K_f} f_k  \prod_{k \in K_1 \cap K_g} g_k, \\
F_2 & = \prod_{k \in K_2 \cap K_f} f_k  \prod_{k \in K_2 \cap K_g} g_k,
\end{align*}

\noindent so that $F=f_2F_1F_2$. By Cauchy's inequality, we have
\begin{equation}
\label{minor-arcs-bounds}
\int_\mathfrak{m} F \ll \sup_{\mathfrak{m}} |f_2| \left( \int_0^1 |F_1|^2 \right)^{1/2} \left( \int_0^1 |F_2|^2 \right)^{1/2}. \\
\end{equation}

\subsection{Weyl's inequality}

To estimate the first term in the bound \eqref{minor-arcs-bounds}, we need pointwise bounds \cite{vaughan_hardy-littlewood_1997} for $f_2$ on the minor arcs. 
\begin{lem}[Weyl's inequality]
\label{lem:weyl-inequality}
Let $\alpha \in \R$. Suppose $(a,q)=1$ and $\alpha \in \R$ such that $|\alpha -a/q| < q^{-2}$. For $k \geqslant 1$ and $\alpha_1, \ldots, \alpha_k \in \C$, let $\phi$ be the polynomial 
\begin{equation}
\phi(x) = \alpha x^k + \alpha_1 x^{k-1} + \cdots + \alpha_{k-1} x + \alpha_k.
\end{equation}

Then we have, for all $Q \geqslant 1$,
\begin{equation}
\sum_{x=1}^Q e(\phi(x)) \ll Q^{1+\varepsilon} \left( q^{-1} + Q^{-1} + qQ^{-k} \right)^{2^{1-k}}.
\end{equation}
\end{lem}

\begin{lem}
For all $\alpha \in \mathfrak{m}$, 
\begin{equation}
f_2(\alpha) \ll f_2(0) n^{-\tau/2+\varepsilon}.
\end{equation}
\end{lem}

\proof This is Weyl's inequality applied to $\phi = \alpha x^2$. We use the fact that in this case $k=2$, $Q\asymp n^{1/2}$ and $n^\tau \leqslant q \leqslant n^{1-\tau}$ by Dirichlet's lemma on Diophantine approximation. We get the result by adding back $f_2(0) \asymp n^{1/2}$. \qed

\subsection{Mean value theorems}

Recall now the mean value theorems of \cite{vaughan_new_1989, vaughan_further_1995} and \cite{wooley_large_1992}, as well as Ford's algorithm \cite{ford_representation_1995}.

\begin{lem}[Vaughan-Wooley]
There is a $\gamma > 0$ so that, for every $k \geqslant 3$ and $s \geqslant 1$, there is a computable $\lambda(k, s) > 0$ such that 
\begin{equation}
\label{mvt}
\int_0^1 |g_k(\alpha)|^{2s} \mathrm{d}\alpha \ll n^{\lambda(k,s)/k + \varepsilon}.
\end{equation}
\end{lem}

Explicit tables of exponents for this mean square theorem are algorithmically computable by the three methods described in \cite[Theorem 4.1]{vaughan_new_1989}, \cite[Lemma 2.3]{vaughan_new_1989-1} and \cite[Lemma 3.2]{wooley_large_1992} optimized in the way described by Wooley \cite{wooley_large_1992}. The implementation of the algorithm in Sage providing all the values used in this paper is available on the authors' webpages. See Table \ref{tablevalues} below for some output of this algorithm, in particular containing the tables of Ford whose algorithm is undisclosed.

\begin{center}
\begin{figure}[h]
\label{tablevalues}
\begin{tabular}{c|c||c|c}
$k$ & $\lambda(k, k)$ & $k$ & $\lambda(k,k)$ \\
\hline
4& 4.60572553279363 &    40& 50.9338839916435 \\
     6& 7.31830866162191   &  42& 53.4919522856964 \\
     8& 9.92905727118400    & 44& 56.0499163246911 \\
    10& 12.5085676596728   &  46& 58.6077897648850 \\
    12& 15.0810335354744   &  48& 61.1655839817793 \\
    14& 17.6492420253841 &    50& 63.7233085263161 \\
    16& 20.2147016775680  &   52& 66.2809714759776 \\
    18& 22.7782942010074   &  54& 68.8385797079435 \\
    20& 25.3405652008671   &  56& 71.3961391137431 \\
    22& 27.9018686743506  &   58& 73.9536547694960 \\
    24& 30.4624435937399  &   60& 76.5111310720912 \\
    26& 33.0224567697859  &   62& 79.0685718489890 \\
    28& 35.5820280054141  &   64& 81.6259804474121 \\
    30& 38.1412454741396   &  66& 84.1833598073007 \\
    32& 40.7001754622901   &  68& 86.7407126613713 \\
    34& 43.2588687351309   &  70& 89.2980408848625 \\
    36& 45.8173648117595   &  72& 91.8553469369745 \\
    38& 48.3756949057251   &  74& 94.4126324955738 
\end{tabular}
\caption{Values of exponents in the mean value theorem.}
\end{figure}
\end{center}

Note that $\lambda(k, s)$ is defined only for integer values of $s$ in the above result. However, by Hölder's inequality we have, for any $\theta \in (0, 1)$ and $k, h \geqslant 3$,
\begin{align*}
\int_0^1 |g_k|^{2(h+ \theta)} & \leqslant \left( \int_0^1 |g_k|^{2h} \right)^{1-\theta} \left( \int_0^1 |g_k|^{2(h+1)} \right)^\theta \\
& \ll n^{\frac{1}{k}\left((1-\theta)\lambda(k, h) + \theta \lambda(k, h+1)\right)}, 
\end{align*}

\noindent so that the relation \eqref{mvt} still holds for any real value of $s$ by letting
\begin{equation}
\lambda(k, h + \theta) = (1-\theta)\lambda(k, h) + \theta \lambda(k, h+1), \qquad \theta \in (0, 1).
\end{equation}

In order to add back the value of $g_k(0)$, we use the fact that the set of smooth numbers is full-sized \cite[Lemma 3.4]{ford_representation_1995} in the sense that, for all $\gamma >0$, we have $|\mathcal{A}(n^{1/k}, n^\gamma)| \gg n^{1/k}$, so that $g_k(0)$ is of size $n^{1/k}$ for all $k \geqslant 1$.  We therefore get for $i \in \{1, 2\}$ and complex numbers $a_k$'s such that $\sum_k a_k^{-1} =1$, 

\begin{align*}
\int_0^1 |F_i|^2 & \ll \prod_{k \in K_i} \left( \int_0^1 |g_k|^{2a_k} \right)^{1/a_k} \ll F_i(0)^2 n^{\phi_i} , \qquad i \in \{1, 2\}.
\end{align*}

\noindent for the constant $\phi_i$ given by
\begin{equation}
\phi_i = \sum_{k \in K_i} \frac{\lambda(k, a_k)}{k a_k} - 2 \sum_{k \in K_i} \frac{1}{k}.
\end{equation}

\textit{Remark.} This formalism will be steadily used all along the paper, the choices for $F_1$ and $F_2$ changing from one section to another, always referring to a partition of the set $K$ except some powers already taken care of.  The choice of $s$ entirely determines the optimal $a_k$ by the optimization algorithm described by Ford, and therefore the exponent $\phi_i$. This exponent decreases with $s$, so that depending on the requirements we can choose the least possible $s$ providing the desired bound. We will repeatedly use this approach in the following, providing the chosen values for $s$ and each specific sets $K_1$ and $K_2$, after making the requirement explicit. We use Ford's heuristics \cite{ford_representation_1996}, claiming that a choice close to the optimal is to take the $a_k's$ such that $a_{k_i} k_j = a_{k_j} k_i$ for all $i, j$ parametrizing $K_1$ or $K_2$.

\subsection{Mixed mean value theorems}
\label{subsec:mvt-mixed}

The problem is that the above mean value theorem only holds in the case of smooth functions. However, the presence of non-smooth functions $f_k$'s are necessary for applying iterative methods. We extend the above mean value theorems by allowing some non-smooth $f_h$. Let $h \in K_f$ and $k_1, \ldots, k_r \in K_g$. \cite{ford_representation_1996} provides an algorithm whose output is a $\phi$  such that
\begin{equation}
\int_0^1 |f_h g_{k_1} \cdots g_{k_r}|^2 \ll n^{\phi},
\end{equation}

\noindent with
\begin{equation}
\phi = \sum_{i=1}^r \frac{x_i}{k_i} \nu(h, k_i, 1/x_i),
\end{equation}

\noindent where the $\nu(h, k_i, x)$ and the optimal values of the convex coefficients $(x_i)_i$ are algorithmically computable. As this algorithm is essential in our paper, we have implemented Ford's algorithm and the code used in this article is provided on the authors'webpages. It can be used for any set of powers $(k_i)_i$.

\subsection{Treatment of minor arcs}
\label{subsec:minor}

All the tools are now at hand to bound the quantities appearing in \eqref{minor-arcs-bounds}. Recall that $\mathfrak{M} = \mathfrak{M}(n^\tau)$ and $\mathfrak{m} = \mathfrak{m}(n^\tau)$ and $\tau$ has to be determined for the circle method to apply while keeping the least possible $s$. The strategy is to use the Cauchy-Schwarz inequality, the Weyl inequality and  the above mean value theorems, writing

\begin{align*}
\int_\mathfrak{m} F & \ll \sup_{\mathfrak{m}} |f_2| \left( \int_0^1 |F_1|^2 \right)^{1/2} \left( \int_0^1 |F_2|^2 \right)^{1/2} \\
& \ll f_2(0) n^{-\tau/2+\varepsilon} \left( F_1(0)^2 n^{\phi_1} \right)^{1/2} \left(  F_2(0)^2n^{\phi_2} \right)^{1/2}.
\end{align*}

\textit{Remarks.} We can shed some light on the choices made for $\tau$ and $s$ :
\begin{itemize}
\item[(i)] Assume for the sake of symmetry and for this heuristic remark that $\phi_1 \simeq \phi_2 \simeq \phi$.  The above bound exponent in $n$ has to be less than $-1$ for the minor arcs contribution to be negligible in from of the main term estimated in Theorem \ref{thm:even-ascending-powers}. It is therefore necessary to have $\phi - \tau/2 < -1$.  In particular, the bound improves for larger $\tau$ and larger $s$ (since increasing $s$ lowers the value of $\phi$).
\item[(ii)] Larger $\tau$ yields better bounds on $s$. Since $\tau$ cannot be larger than $1/2$, this implies that we need $\phi < -3/4$, what provides a lower bound on $s$ with the chosen method, namely around $s=72$.
\item[(iii)] The final choice is a trade-off between the quality of the bound on minor arcs and the quality of approximation on major arcs in Section \ref{sec:approximation}, in order to get the least possible $s$. We get the following table when optimizing the choice of partitions of the form $K_1 = \{6, 8, \ldots, n\}$ and $K_2 = \{4, n+2, n+4, \ldots, 2s\}$  :

\begin{center}
\begin{figure}[h]
\begin{tabular}{c|c|c}
$\tau$ & $2s$ & $n$ \\
\hline\\[-.8em]
0.385 & 278 & 54  \\
0.386 & 276 & 54  \\
0.387 & 274 & 54 \\
0.388 & 272 & 52 \\
0.389 & 270 & 52 \\
0.390 & 268 & 52 \\
0.392 & 266 & 52 \\
0.393 & 264 & 52 \\
0.394 & 262 & 52 
\end{tabular}
\caption{Numerical relations between $\tau$ and $s$ for bounding minor arcs.}
\end{figure}
\end{center}
\end{itemize}

From now on, let $\tau = 0.3935$ and $2s=266$. We get the following :

\begin{lem}
There is $\gamma >0$ such that
\begin{equation}
\int_{\mathfrak{m}} F \ll F(0)n^{-1-\delta}.
\end{equation}
\end{lem}

\proof The previous algorithm yields the values
\begin{align*}
\phi_1 & = -0.806 \quad \text{with} \quad K_1 = \{6, 8, \ldots, 52\}, \\
\phi_2 & = -0.801 \quad \text{with} \quad  K_2 = \{4,54, 56, \ldots, 266\}.
\end{align*}

By Cauchy-Schwarz inequality and the previous bounds we get, for $\delta = 0.0007$,
\begin{align*}
\int_\mathfrak{m} F & \ll \sup_{\mathfrak{m}} |f_2| \left( \int_0^1 |F_1|^2 \right)^{1/2} \left( \int_0^1 |F_2|^2 \right)^{1/2} \\
& \ll f_2(0) n^{-\tau/2+\varepsilon} \left( F_1(0)^2 n^{\phi_1} \right)^{1/2} \left(  F_2(0)^2n^{\phi_2} \right)^{1/2} \\
& \ll F(0)n^{-1-\delta}. \qed
\end{align*}

We therefore deduce the following reduction of the problem : 

\begin{lem}
There is $\delta > 0$ such that
\begin{equation}
\int_0^1 F(\alpha)e(-n\alpha)\mathrm{d}\alpha = \int_{\mathfrak{M}}   F(\alpha)e(-n\alpha)\mathrm{d}\alpha + O\left(F(0)n^{-1-\delta}\right).
\end{equation}
\end{lem}

\section{Approximation on major arcs}
\label{sec:approximation}

\subsection{Approximated versions}

We introduce the approximated versions of the generating functions $f_k$ and $g_k$. Introduce $e_q(x) = e(x/q)$. In this whole section, all the $\alpha$ considered are in $\mathfrak{M}$. Note $\rho$ the Dickman's function introduced in \cite[Section 4]{ford_representation_1996}. For $k \geqslant 1$, let $(a,q)=1$ and $|\beta| < 1/2$ such that $\alpha = \frac{a}{q} + \beta$. Define
\begin{align*}
S_k(q,a) & = \sum_{m=1}^q e_q\left( a m^k \right), \\
w_k(\beta) & = \sum_{m \leqslant n} k^{-1} m^{1/k-1} e(\beta m), \qquad \text{for } \ k \leqslant 4, \\
w_k(\beta) & = \sum_{n^{\gamma k} < m \leqslant n} k^{-1} m^{1/k-1} \rho\left( \frac{\log m}{\gamma k \log n} \right) e(\beta m), \qquad \text{for } \ k \geqslant 5, \\
W_k(\alpha) & = q^{-1} S_k(q,a) w_k(\beta), \\
\Delta_k(\alpha) & = f_k(\alpha) - W_k(\alpha).
\end{align*}

We want to replace the $f_k$ and $g_k$ by $W_k$ for each $k \geqslant 1$, up to an error term of size $F(0) n^{-1-\delta}$. To this end, recall the following two pointwise bounds \cite[Lemmas 4.1 and 4.2]{ford_representation_1995}.

\begin{lem}
\label{lem:bound-W}
Let $k \geqslant 1$. For all $(a,q)=1$ and $|\beta|<1/2$ we have, for $\alpha = \frac{a}{q} + \beta$, 
\begin{equation}
W_k(\alpha) \ll \left( \frac{n}{q} \right)^{1/k} (1+n |\beta|)^{-1}.
\end{equation}
\end{lem}

\begin{lem}
\label{lem:error-delta}
Let $k \geqslant 1$. For all $(a,q)=1$ and $|\beta| < 1/2$ we have, for $\alpha = \frac{a}{q} + \beta$, 
\begin{equation}
\Delta_k(\alpha) \ll q^{1/2+\varepsilon} (1+n|\beta|)^{1/2}.
\end{equation}
\end{lem}

\subsection{Replacing $f_2$}

In order to replace $f_2$ by its approximated version $W_2$, we recall the bound on the approximation error in Lemma \ref{lem:error-delta}, yielding
\begin{equation}
\Delta_2 \ll n^{\tau/2+\varepsilon} \ll f_2(0) n^{(\tau-1)/2 + \varepsilon}.
\end{equation}

Therefore, by the mean value bounds stated in \eqref{mvt} with the same $F_1$ and $F_2$ as in the previous section, taking the same $\delta$ as before, we have
\begin{align*}
\int_{\mathfrak{M}} |\Delta_2 f_4f_6g_8g_{10}F_2| & \ll f_2(0) n^{\frac{\tau - 1}{2} + \varepsilon} \left( \int_0^1 |F_1|^2 \right)^{1/2} \left( \int_0^1 |F_2|^2 \right)^{1/2} \\
& \ll F(0)n^{(\tau-1)/2 + \phi_1/2 + \phi_2/2 + \varepsilon} \\
& \ll F(0)n^{-1-\delta}.
\end{align*}

The problem is therefore reduced as follow :

\begin{lem}
There is $\delta > 0$ such that
\begin{equation}
\int_{\mathfrak{M}}  F(\alpha)e(-n\alpha)\mathrm{d}\alpha = \int_{\mathfrak{M}}   W_2 f_4 \cdots g_{2s} e(-n\alpha)\mathrm{d}\alpha + O\left(F(0)n^{-1-\delta}\right).
\end{equation}
\end{lem}

\subsection{Replacing $f_4$}
\label{f4}

Recall the pruning lemma  \cite{brudern_problem_1988}, central in the approximation process.

\begin{lem}[Br\"{u}dern, pruning lemma]
Let $X \leqslant n$. For $1 \leqslant a \leqslant q \leqslant X$ with $(a, q)=1$, let $\mathfrak{M}(X ; q, a)$ be an interval in $((a/q) - 1/2, (a/q) + 1/2)$ and assume the $\mathfrak{M}(X ; q, a)$ are pairwise disjoint. Let $G : \mathfrak{M}(X) \to \C$ be a function such that
\begin{equation}
G(\alpha) \ll \frac{n}{q} (1+n\beta)^{-1}, \qquad \text{for } \ \alpha \in \mathfrak{M}(X ; q, a).
\end{equation} 

Let $\Psi : \R \to [0, +\infty)$ be a function with Fourier expansion of the form
\begin{equation}
\Psi(\alpha) = \sum_{|h| \ll n^\eta} \psi_h e(\alpha h),
\end{equation}

\noindent for a certain $\eta > 0$, and such that 
\begin{equation}
\psi_0 = \int_0^1 \Psi(\alpha) \mathrm{d}\alpha \ll X^{-1}\Psi(0).
\end{equation}

Then
\begin{equation}
\int_{\mathfrak{M}(X)} G(\alpha) \Psi(\alpha)\mathrm{d}\alpha \ll n^\varepsilon \Psi(0).
\end{equation}
\end{lem} 

The strategy is to apply this pruning lemma and the Cauchy inequality to write
\begin{align*}
\int_{\mathfrak{M}} |W_2 \Delta_4 F_1 F_2| & \ll \sup_{\mathfrak{M}} |\Delta_4| \left( \int_{\mathfrak{M}} |W_2^2 F_1| \right)^{1/2} \left( \int_0^1 |F_2^2| \right)^{1/2} \\
& \ll n^{\tau/2+\varepsilon} (n^{\varepsilon} F_1(0)^2)^{1/2} \left( F_2(0)^2 n^{\phi_2}\right)^{1/2} .
\end{align*}

\textit{Remarks.} This strategy motivates some heuristic comments :

\begin{itemize}
\item[(i)] To satisfy the assumption of the pruning lemma, we need to choose $K_1$ such that $\phi_1 < -\tau$. Moreover, for the bound to be sufficient for our purposes, we need the final exponent to satisfy $\tau + \phi_2 < -1/2$. These conditions add up to an optimization problem that ultimately justifies the choices made for $\tau$, $s$ and the sets $K_1$, $K_2$. 
\item[(ii)] Unlike the minor arcs situation, larger $\tau$ yield worse bounds since the requirement on $\phi_1$ would be stronger. We chose $\phi_1$ to be as close as possible to $-\tau$ ; and $\phi_2$ is chosen so that the second bound is satisfied.
\item[(iii)] We get the following table when considering sets of the form $K_1 = \{22, \ldots, n\}$ and $K_2 = \{6, \ldots, 20,  n+2, n+4, \ldots, 2s\}$ :

\begin{center}
\begin{figure}[h]
\begin{tabular}{c|c|c}
$\tau$ & $2s$ & $n$ \\
\hline\\[-.8em]
0.390 & 256 & 50 \\
0.391 & 258 & 50 \\
0.392 & 260 & 50 \\
0.393 & 264 & 50 \\
0.394 & 266 & 50 \\
0.395 & 268 & 50 \\
0.396 & 282 & 52
\end{tabular}
\caption{Numerical relations between $\tau$ and $s$ for replacing $f_4$.}
\end{figure}
\end{center}

In particular Tables 2 and 3 shows that the choice of $\tau = 0.3935$ allows the value of $2s=266$. 
\end{itemize}

Let
\begin{align*}
K_1 & = \{22, \ldots, 50\},  \\
K_2 & = \{6, \ldots, 20, 52, \dots, 268\}.
\end{align*}

In order to replace $f_4$ by $W_4$, we want to use  this pruning lemma with $G=W_2^2$ and $\Psi = |F_1|^2 =  |g_{22}  \cdots g_{50}|^2$. The Vaughan iterative method yields $\phi_1 < - \tau$, \textit{i.e.}
\begin{equation}
\label{re}
\int_0^1 \Psi \ll \Psi(0) n^{-\tau}, 
\end{equation}

\noindent  so that the pruning lemma applies and gives
\begin{equation}
\int_{\mathfrak{M}} |W_2 F_1|^2 \ll n^\varepsilon F_1 (0)^2.
\end{equation}

Moreover, with $F_2 = f_6 g_8  \cdots g_{20} g_{52} \cdots g_{268}$, the mean value theorem yields the required value in order to have a small enough bound below, namely 
\begin{equation}
\int_0^1 |F_2|^2 \ll F_2(0)^2 n^{-0.895}.
\end{equation}

We therefore can conclude by Cauchy-Schwarz inequality and with the bound for $\Delta_4$ given in Lemma \ref{lem:error-delta}, with $\delta = 0.0007$,
\begin{align*}
\int_{\mathfrak{M}} |W_2 \Delta_4 F_1 F_2| & \ll \sup_{\mathfrak{M}} |\Delta_4| \left( \int_0^1 |W_2^2 F_1| \right)^{1/2} \left( \int_0^1 |F_2^2| \right)^{1/2} \\
& \ll n^{\tau/2+\varepsilon} (n^{\varepsilon} F_1(0)^2)^{1/2} \left( F_2(0)^2 n^{-0.84}\right)^{1/2}  \\
& \ll F(0) n^{-1-\delta}.
\end{align*}

\begin{lem}
There is $\delta > 0$ such that
\begin{equation*}
\int_{\mathfrak{M}}   W_2 f_4 \cdots f_{2s} e(-n\alpha)\mathrm{d}\alpha = \int_{\mathfrak{M}}   W_2 W_4 f_6 \cdots f_{2s} e(-n\alpha)\mathrm{d}\alpha + O\left(F(0) n^{-1-\delta}\right).
\end{equation*}
\end{lem}

\section{Pruning major arcs}
\label{sec:pruning}

In  order to be able to replace $f_6$ by $W_6$, we have to prune the major arcs.

\subsection{Pruning to $n^{\kappa}$}

Introduce the new major arcs $\mathfrak{M}_1 = \mathfrak{M}(n^{\kappa})$ for a $\kappa < \tau$ and the associated relative minor arcs $\mathfrak{m}_1 = \mathfrak{M} \backslash \mathfrak{M}_1$. We will take $\kappa = 1/4$ for reasons that will appear in Section \ref{subsec:f6}. We have, by Lemma \ref{lem:bound-W},
\begin{equation}
|W_2W_4|^2 \ll \left(\frac{n}{q}\right)^{3/2} (1+n\beta)^{-4}.
\end{equation}

Examine more precisely the elements in $\mathfrak{m}_1$. For any $0 < \theta < \kappa$, we have either $q \geqslant n^{ \kappa-\theta}$ (case I) or $\beta \geqslant n^{\theta-1}$ (case II). We therefore have one of the two bounds
\begin{align*}
(\text{case I}) \qquad & |W_2W_4|^2 \ll \frac{n^{3/2}}{q} (1+|\beta| n )^{-1} n^{(\theta - \kappa)/2}, \\
(\text{case II}) \qquad & |W_2W_4|^2 \ll \frac{n^{3/2}}{q} (1+|\beta|n )^{-1} n^{-3\theta}.
\end{align*}

Taking $\theta = \kappa/7$  in order to make both bounds match we get that, on $\mathfrak{m}_1$,
\begin{equation}
|W_2W_4|^2 \ll \frac{n}{q} (1+|\beta|n )^{-1} n^{1/2-3\theta}.
\end{equation}

In particular the pruning lemma can be applied to $G = n^{3\theta - 1/2}|W_2W_4|^2$. Moreover, redefining $F_1 = g_{22} \cdots g_{50}$, we can apply the pruning lemma since by the mean value theorems we have 
\begin{equation}
\int_0^1 |F_1|^2 \ll F_1(0)^2 n^{-0.3958}.
\end{equation}

Moreover, applying the mean value algorithm we get, with $F_2 = f_6 g_8 \cdots g_{20} g_{52} \cdots g_{266}$, 
\begin{equation}
\int_0^1 |F_2|^2 \ll F_2(0)^2 n^{-0.895}.
\end{equation}

Altogether, the contribution of the minor arcs are shown to be negligible, since by Cauchy-Schwarz inequality we get, letting $\delta = 0.001$, 
\begin{align*}
\int_{\mathfrak{m}_1} |W_2W_4f_6 g_8 \cdots g_{266}|^2 & 
\ll n^{1/4-3\theta/2} \left( \int_{\mathfrak{m}_1} n^{3\theta - 1/2}|		W_2W_4F_1|^2 \right)^{1/2} \left( \int_{0}^1 |F_2|^2 \right)^{1/2} \\
& \ll F(0) n^{1/4-3\theta/2 - 1/2 - 1/4+ \phi_2/2} \\
& \ll F(0)n^{-1-\delta}.
\end{align*}

In particular it is sufficient to study the integral formulation of the problem with the integral restricted to the new major arcs $\mathfrak{M}_1$ :

\begin{lem}
There is $\delta > 0$ such that
\begin{equation*}
 \int_{\mathfrak{M}}   W_2 W_4 f_6 g_8 \cdots g_{2s} e(-n\alpha)\mathrm{d}\alpha =  \int_{\mathfrak{M}_1}   W_2 W_4 f_6 g_8 \cdots g_{2s} e(-n\alpha)\mathrm{d}\alpha + O\left(F(0) n^{-1-\delta}\right).
\end{equation*}
\end{lem}

\subsection{Replacing $f_6$}
\label{subsec:f6}

Now that we are reduced to major arcs of size $n^{\kappa}$ for $\kappa = 1/4$, it is possible to replace $f_6$. In order to efficiently apply the pruning lemma, we distribute $W_k^k$ in Hölder's inequality. We get
\begin{align*}
& \int_{\mathfrak{M}_1} |W_2W_4\Delta_6f_8 \cdots f_{2s}| \\
& \qquad \ll  \sup_{\mathfrak{M}_1} |\Delta_6| \left( \int_{\mathfrak{M}_1} |W_2^2F_1^2| \right)^{1/2} \left( \int_{\mathfrak{M}_1} |W_4^4F_2^2| \right)^{1/4} \left( \int_0^1 |F_2^2| \right)^{1/4} \\
& \qquad \ll F(0)n^{\kappa/2 -1/2 - 1/4 - 1/6 + \phi_2/4}.
\end{align*}

\textit{Remark.} It is because of this lemma that the major arcs have been pruned above to $n^{\kappa}$ with $\kappa = 1/4$. Indeed, assuming $s$ is large enough for the pruning lemma's hypotheses to be satisfied, the bound can reach $F(0)n^{-1-\delta}$ if and only if $\kappa < 1/4$. Moreover, the closer we are of $1/3$, the closer $\phi$ has to be to one, \textit{i.e.} the larger $s$ has to be. This justifies the choice of $\kappa$ slightly away from $1/3$.

The above bound is negligible as soon as $\phi_2 < -5/6$ and the pruning lemma is indeed applicable as soon as $\phi_1, \phi_2 < - \kappa$. This is the case with 
\begin{align*}
K_1 & = \{64, \ldots, 108\}, \\
K_2 & = \{8, \ldots, 62, 110, \ldots, 266\}.
\end{align*}

With this choice we get
\begin{align*}
\int_0^1 |F_1|^2 & \ll F_1(0)^2 n^{-0.255}, \\
\int_0^1 |F_2|^2 & \ll F_2(0)^2 n^{-0.897}.
\end{align*}

Finally, we get, with $\delta = 0.01$,
\begin{align*}
\int_{\mathfrak{M}_1} |W_2W_4\Delta_6g_8 \cdots g_{266}| \ll F(0)n^{-1 - \delta}.
\end{align*}

As such, this contribution is negligible, and the problem further reduces as follows.

\begin{lem}
There is $\delta > 0$ such that
\begin{equation*}
 \int_{\mathfrak{M}_1}   W_2 W_4 f_6 g_8 \cdots g_{2s} e(-n\alpha)\mathrm{d}\alpha =  \int_{\mathfrak{M}_1}   W_2 W_4 W_6 g_8 \cdots g_{2s} e(-n\alpha)\mathrm{d}\alpha + O\left(F(0) n^{-1-\delta}\right).
\end{equation*}
\end{lem}

\subsection{Pruning to $\log^A n$}

Let $Y = \log^A n$ for a certain $A >0$. Introduce the logarithmically-pruned major arcs $\mathfrak{M}_2 = \mathfrak{M}(Y)$, and $\mathfrak{m}_2$ the associated related minor arcs that is to say $\mathfrak{m}_2 = \mathfrak{M}_1 \backslash \mathfrak{M}_2$. As for the previous section, we begin by showing that the contribution of these new minor arcs is negligible. On $\mathfrak{m}_2$, by Lemma \ref{lem:bound-W}, we have $W_2^2, W_4^4, W_6^6 \ll \frac{n}{q} (1+n\beta)^{-1}$. Choosing
\begin{align*}
  K_1 &= \{8, 204, 206, \ldots, 266 \},\\
  K_2 &= \{10,12,40, 42, \ldots, 48 \},\\
  K_3 &= \{14,16,\ldots,38,50,52,\ldots,202 \},
\end{align*}
the pruning lemma applies with the bounds 
\begin{align*}
\int_0^1 |F_1|^2 & \ll F_1(0)^2  n^{-0.252}, \\
\int_0^1 |F_2|^2 & \ll F_2(0)^2  n^{-0.286}, \\
\int_0^1 |F_3|^2 & \ll F_3(0)^2 n^{-0.823}.
\end{align*}
Together with
\[ 
  \int_0^1 |F_2|^6 |F_3|^8 \ll F_2(0)^6 F_3(0)^8 n^{-1.0002},
\]

and $\delta = 0.00001$, we get

\begin{align*}
& \int_{\mathfrak{m}_2} |W_2W_4W_6F_1F_2F_3| \\
&  \ll \left( \int_{\mathfrak{m}_2} |W_2F_1|^2 \right)^{1/2}  \left( \int_{\mathfrak{m}_2} |W_4^4F_2^2| \right)^{1/4} \left( \int_{\mathfrak{m}_2} |W_6^6F_3^2| \right)^{1/6} \left( \int_{\mathfrak{m}_2} |F_2^6F_3^8| \right)^{1/12} \\
& \ll F(0) n^{-1-\delta}.
\end{align*}

We therefore reduced the problem to estimating the integral on the new arc $\mathfrak{M}_2$.

\begin{lem}
There is $\delta > 0$ such that
\begin{equation*}
 \int_{\mathfrak{M}_1}   W_2 W_4 W_6 \cdots g_{2s} e(-n\alpha)\mathrm{d}\alpha =  \int_{\mathfrak{M}_2}   W_2 W_4 W_6 \cdots g_{2s} e(-n\alpha)\mathrm{d}\alpha + O\left(F(0) n^{-1-\delta}\right).
\end{equation*}
\end{lem}

\subsection{Pruning to $\log^{1/4} n$}

It is  now necessary to prune again these major arcs to reach logarithmically sized arcs, so that all the remaining $g_k$'s will be directly approachable without efforts by $W_k$. Introduce $Z = \log^{1/4} n$ and the new major and minor arcs $\mathfrak{M}_3 = \mathfrak{M}(Z)$ and $\mathfrak{m}_3 = \mathfrak{M}_2 \backslash \mathfrak{M}_3$ the associated relative minor arcs. On $\mathfrak{m}_3$, the Vaughan-Wooley bound \cite{vaughan_further_1995} implies the following.
\begin{lem}[Vaughan-Wooley]
We have, for all $\alpha \in \mathfrak{m}_3$, 
\begin{equation}
g_k(\alpha) \ll \left( \frac{n}{q} \right)^{1/k} q^\varepsilon (1+n\beta)^{-1/k}.
\end{equation}
\end{lem}

In particular we get, by Lemma \ref{lem:bound-W},
\begin{equation}
|W_2W_4W_6g_8 \cdots g_{2s}| \ll F(0) q^{-\mu} (1+n\beta)^{-\eta},
\end{equation}

\noindent where 
\begin{align*}
\omega & = \frac{1}{2} + \cdots + \frac{1}{2s} - \varepsilon,  \\
\eta & = 3 + \frac{1}{8} + \cdots + \frac{1}{2s}.
\end{align*}

Now, integrating over the minor arcs $\mathfrak{m}_3$ we get
\begin{align*}
\int_{\mathfrak{m}_3} |W_2W_4W_6g_8 \cdots g_{2s}| & \ll F(0)\left( \sum_{q \leqslant Z} q^{1-\omega} \int_{Z/nq}^\infty  \frac{\mathrm{d}\beta}{(n\beta)^\eta} + \sum_{Z < q \leqslant Y}  q^{1-\omega} \int_0^{Y/nq}  \frac{\mathrm{d}\beta}{(1+n\beta)^\eta} \right) \\
& \ll F(0) \left( n^{-1} Z^{1-\eta} \sum_{q \leqslant Z} q^{\eta-\omega} + n^{-1} \sum_{Z < q \leqslant Y} q^{1-\omega} \right) \\
& \ll F(0) n^{-1} Z^{2-\omega}.
\end{align*}

Altogether, we get for $\rho = (2-\omega)/4 >0$,
\begin{equation}
\int_{\mathfrak{m}_3} |W_2W_4W_6g_8 \cdots g_{2s}| \ll F(0)n^{-1} \log^{-\rho} n,
\end{equation}

\noindent which is negligible in front of the expected main term estimated in Theorem \ref{thm:even-ascending-powers}. In particular, it is enough to concentrate on the integral whose integration domain is restricted to $\mathfrak{M}_3$ :

\begin{lem}
There is $\alpha > 0$ such that
\begin{equation*}
 \int_{\mathfrak{M}_2}   W_2 W_4 W_6 \cdots g_{2s} e(-n\alpha)\mathrm{d}\alpha =  \int_{\mathfrak{M}_3}   W_2 W_4 W_6 \cdots g_{2s} e(-n\alpha)\mathrm{d}\alpha + O\left(F(0)n^{-1} \log^{-rho} n\right).
\end{equation*}
\end{lem}

\subsection{Replacing the remaining $g_k$'s}

On these log-sized major arcs, it is straightforward to replace $g_k$ by $W_k$ up to an error term as soon as $k \geqslant 5$. Indeed, we have \cite[Equation (4.21)]{ford_representation_1996} for all $k \geqslant 5$, 
\begin{equation}
\Delta_k(\alpha) \ll n^{1/k} \log^{-3/4} n.
\end{equation}

Using the fact that the major arcs $\mathfrak{M}_3$ are of length $Z^2n^{-1}$, we get the bound
\begin{equation}
\int_{\mathfrak{M}_3} |\Delta_k| \ll n^{1/k-1} \log^{-1/4} n \ll g_k(0) n^{-1} \log^{-1/4} n.
\end{equation}

With the bounds given in Lemmas \ref{lem:weyl-inequality} and \ref{lem:bound-W}, we deduce for any $r \geqslant 4$, 
\begin{align*}
\int_{\mathfrak{M}_3} W_2 \cdots W_{2r-2} \Delta_{2r} g_{2r+2} \cdots g_{2s} \ll F(0) n^{-1} \log^{-1/4} n.
\end{align*}

We therefore deduce, adding finitely many such approximation error terms :

\begin{lem}
We have
\begin{align*}
& \int_{\mathfrak{M}_3} W_2 \cdots W_6 g_8 \cdots g_{2s}  e(-n\alpha) \mathrm{d}\alpha \\
& \qquad \qquad  = \int_{\mathfrak{M}_3} W_2 \cdots W_{2s} e(-n\alpha) \mathrm{d}\alpha + O\left(F(0) n^{-1} \log^{-1/4} n\right).
\end{align*}
\end{lem}

\textit{Remark.} Unlike the previous works of Thanigasalam \cite{thanigasalam_additive_1968} or Ford \cite{ford_representation_1995}, all the $g_k's$ have been replaced by their approximated versions $W_k$, greatly improving the treatment of the singular quantities below.

\subsection{Completing the arcs}

The last step before studying the expression yielding the main term is to replace the integration over each $\mathfrak{M}(Z; q, a)$ above by the integral over the whole circle. This is doable since $W_k$ is very small outside $\mathfrak{M}(Z; q, a)$. We follow the method of \cite[(4.33)]{ford_representation_1996} and get that the complementary part is
\begin{align*}
 \sum_{q \leqslant Z} \sum_{(a, q) = 1} \int_{[0,1] \backslash \mathfrak{M}(Z; q, a)} |W_2 \cdots W_{2s}(\alpha ; q, a)|& \ll F(0) \sum_{q \leqslant Z} q^{1-\omega} \int_{Z/qn}^{+\infty} \frac{\mathrm{d}\beta}{(n\beta)^\eta} \\
& \ll F(0) n^{-1} \log^{-\rho} n.
\end{align*}

Altogether, we proved the following:

\begin{lem}
We have, for $Z = \log^{1/4} n$,
\begin{align*}
& \int_{\mathfrak{M}_3} W_2 \cdots W_{2s} e(-n\alpha) \mathrm{d}\alpha  \\
& \qquad = \mathfrak{S}(n, Z) I(n) + O\left(F(0) n^{-1} \log^{-1/4} n\right),
\end{align*}

\noindent where
\begin{align*}
\mathfrak{S}(n, Z) & = \sum_{q \leqslant Z} A(n, q), \\
A(n, q) & = \sum_{(a, q)=1} q^{-s} S_2 \cdots S_{2s} (q, a) e_q(-an), \\
I(n) & = \frac{1}{2^s s!} \sum_{(\star)} m_2^{-1/2} m_4^{-3/4} m_6^{-5/6} \prod_{k=4}^s \rho\left( \frac{\log m_k}{k \gamma \log n} \right) m_k^{1/k-1},
\end{align*}

\noindent and the conditions $(\star)$ of the above sum are given by
\begin{align*}
& \frac{n}{2^k} < m_k \leqslant n, \qquad \text{for } k = 2, 4, 6, \\
& n^{\gamma k} < m_k \leqslant n, \qquad \text{for } k \geqslant 8, \\
& n = m_2+ m_4 + \cdots + m_{2s}.
\end{align*}
\end{lem}

\section{Singular series and integral}
\label{sec:singular}

The treatment of the singular series and integral is simplified by the replacement of all the $f_k$ and $g_k$ by their approximated version $W_k$. The argument is analogous to \cite{ford_representation_1996} and we briefly include the details for completeness.

Since the $\displaystyle \frac{\log m_k}{k \gamma \log n}$ are uniformly bounded, the explicit definition of $I(n)$ yields
\begin{equation}
I(n) \gg F(0) n^{-1}. 
\end{equation}

Now it remains to prove that $\mathfrak{S}(n, Z) \gg 1$. First of all, it is straightforward to see that this is a convergent series since $S_k(q,a) \ll q^{1-1/k}$ by \cite[Theorem 4.2]{vaughan_hardy-littlewood_1997}. This henceforth yields
\begin{equation}
|\mathfrak{S}(n, Z) - \mathfrak{S}(n)| \leqslant \sum_{q > Z} |A(n, q)| \ll \sum_{q > Z} q^{1-\omega} \ll \log^{-0.18} n.
\end{equation}

By \cite[Lemma 2.11]{vaughan_hardy-littlewood_1997} the function $A(n,q)$ is multiplicative in $q$. Therefore the study is reduced to the associated local factors, namely we can write
\begin{equation}
\mathfrak{S}(n) = \sum_{q=1}^\infty A(n, q) = \prod_p \chi_p, \qquad \text{where} \quad \chi_p = \sum_{h=0}^\infty A\left(n, p^h\right).
\end{equation}

The same bound as above therefore yields
\begin{equation}
|\chi_p - 1| \leqslant \sum_{h=1}^\infty \left|A\left(n, p^h\right)\right| \ll p^{1-\omega}.
\end{equation}

With the bound $\chi_p \gg p^{-1188}$, consequence of \cite[Lemma 6.4]{ford_representation_1995} with $r=133$ and $\gamma = 9$ in the present case, this leads to $\mathfrak{S}(n) \gg 1$. In particular, $r(n) > 0$ for $n$ large enough, achieving the proof of Theorem \ref{thm:even-ascending-powers}.

\newpage

\bibliographystyle{alpha}
\bibliography{biblio_growing_powers_study}

\end{document}